\documentstyle[amscd,amssymb,12pt]{amsart}
\input xypic \xyoption{curve}
\input epsf

\def\hcorrection#1{\advance\hoffset by #1 }
\def\vcorrection#1{\advance\voffset by #1 }

\vcorrection{-.7in}
\hcorrection{-.5in}

\newcommand{\B}[1]{{\bold#1}} 
\newcommand{\C}[1]{{\mathcal#1}} 

\newcommand{\g}{{\mathfrak g}} 
 
\theoremstyle{plain}
\newtheorem{th}{Theorem}[section]

\newtheorem{prop}{Proposition}[section]

\theoremstyle{definition}
\newtheorem{defin}{Definition}[section]

\theoremstyle{definition}

\theoremstyle{remark}
\newtheorem{rem}{Remark}[section]

\numberwithin{equation}{section}

\begin{document}

\title{From Lie Theory to Deformation Theory\\ and\\Quantization}
\author{Lucian M. Ionescu}
\address{Department of Mathematics, Illinois State University, IL 61790-4520, USA}
\email{lmiones@@ilstu.edu}%

\keywords{Deformation theory, Kuranishi functor, quantization, renormalization}
\subjclass{Primary:14Dxx; Secondary: 46Lxx,53Dxx} 
\date{\today}

\begin{abstract}
{\em Deformation Theory} is a natural generalization of {\em Lie Theory},
from Lie groups and their linearization, Lie algebras,
to differential graded Lie algebras and their higher order deformations, 
quantum groups.

The article focuses on two basic constructions of deformation theory:
the universal solution of Maurer-Cartan Equation (MCE),
which plays the role of the exponential of Lie Theory, 
and its inverse, the Kuranishi functor, as the logarithm.

The deformation functor is the gauge reduction of MCE, corresponding 
to a Hodge decomposition associated to the strong deformation retract data.

The above comparison with {\em Lie Theory} leads to a better understanding 
of {\em Deformation Theory} and its applications, 
e.g. the relation between quantization and Connes-Kreimer renormalization,
quantum doubles and Birkhoff decomposition.
\end{abstract}

\maketitle
\tableofcontents

\section{Introduction}
{\em Lie Theory} (LT) is arguably a model of a well-designed and fundamental theory \cite{VBLT,ELT}.
In contrast, ``Deformation Theory'' (DT) does not have such well-defined ``boundaries'',
i.e. a clearly delimited core of concepts and results at its foundations,
being currently rather a methodology of deriving new structures by deforming ``old'' ones.
This makes some of modern results in quantum mathematical-physics appear as ``independent'' theories,
e.g. in renormalization, while being in fact well ``disguised'' applications 
of deformation theory \cite{IM}.

In this article we aim to compare {\em Lie Theory} and {\em Deformation Theory},
starting with the review of the theory of differential graded Lie algebra structures (DGLA),
as an obvious generalization of Lie algebras.

The further (technical) generalization to L-infinity algebras is straightforward.
On the other hand a comprehensive treatment of this project 
is beyond the author's present capabilities and expertise.

The present investigation is centered on the formal solutions of
Maurer-Cartan equation (MCE), as a formal {\em
Feynman-Taylor power series} \cite{Kon1,L-inf,CFG}.  
It is a representation of an exponential map,  as proved in
\cite{FI} for the special case of the Cartan-Eilenberg complex of $T_{poly}$ and $D_{poly}$.
This suggests that the Kuranishi functor, its inverse,  is the DT-analog of the
logarithm of Lie Theory. 

The main references regarding the universal solutions of the MCE and the Kuranishi
functor, including a background on deformation theory, are \cite{HS,Man}.

The construction of a universal solution of the Maurer-Cartan solution
of a DGLA from \cite{HS} (see also \cite{AI-HS}) is explained in \S\ref{S:dtusol},
in analogy with Picard's method for solving differential (integral) equations,
by using a section of the differential and an iterative procedure (series expansion).
This allows to establish a parallel with Lie Theory,
with emphasis on the role of the {\em contraction} \S\ref{S:hlt},
paving the way towards Rota-Baxter algebras and renormalization \cite{I-RBA}.

As a new result in this article, 
a Hodge decomposition (Definition \ref{D:hd}) is obtained in Theorem \ref{T:splitting},
underlying the strong deformation data (SDR) framework of \cite{HS}
and suggested by a previous observation \cite{AI-HS}. 
In \S\ref{S:hd} it is proved that the contraction of \cite{HS} provides a *-operator, 
together with the associated Dirac operator as a square root of the Laplacian.
This aspect, in relation with Hodge Theory will be developed elsewhere.

The deformation theory point of view of \cite{Man} is adopted in \S\ref{S:dt},
to show that the universal twisting cocycle plays the role of the 
exponential of LT, while its inverse is the Kuranishi map,
the logarithm of LT.
Together with the Kuranishi map, it prompt to an analogy with Lie theory,
relating the infinitesimal and global structures: Theorems \ref{T:log, T:main}.

The present results lead to a few notable consequences.
On the more technical side, 
the connection with generalized complex structures \cite{GCS} is made 
via the interpretation of the above *-operator as a complex structure 
on the corresponding complement of harmonic forms,
and a comparison with the framework characteristic of the $dd^*$-lemma \cite{dd}
(again, to be developed elsewhere).

The second part of the article is concerned with the conceptual implications to
quantization and renormalization, starting with the
general philosophy of deformation quantization \S\ref{S:dorh}: 
to double or to half (i.e. factorize) the algebraic structure,
e.g. Drinfeld's double, Manin triples or Rota-Baxter algebras and Birkhoff decomposition.

The role of the deformation of structure associated to
a contraction, as part of the SDR data, 
is to split the epimorphism and decompose the boundaries $b=\partial c$ (exact cycles)
into source and target components $b^\pm=\partial^\pm c$, 
characteristic of non-abelian cohomology \cite{NACOH},
and a framework characteristic of bialgebra deformation quantization.

The concluding section further clarifies the role of quantization by deformation in general:
it is a {\em categorification} in ``disguise'' (the element $b$ is represented as  
the morphism $b:b^-\to b^+$ \cite{Cat,Cat-bar}), 
anticipated in \cite{I-Q++}, and to be explained in detail elsewhere \cite{I-RBA}.
This also explains the advantages of the Feynman Path Integral approach to quantization,
over the traditional deformation quantization approach.
In particular, the relation between deformation theory and Connes-Kreimer renormalization \cite{CK}
established in this article, will be detailed as a part of the {\em Theory of Rota-Baxter Algebras},
as a generalization of the {\em Theory of Hopf algebras} and {\em Quantum groups},
towards the natural framework of {\em Feynman Processes} \cite{I-PROPF,I-FL}.


\section{Recall on Lie Theory and its subsequent developments}\label{SS:LT}
The main ideas of Lie Theory are summarized following \cite{ELT}, \S\S 3.2, p.82-83
(see also \cite{VBLT}).
The subsequent developments by Cartan, Chen and Kodaira, are summarized from \cite{wiki:FG}.
The connection with Rota-Baxter algebras and Hopf algebra renormalization
leads to the designing principles of path integrals.

\subsection{Overview of Lie Theory}\label{SS:OLT}
The motivation of Sophus Lie's work was to solve differential equations
using groups of transformations. 
His aim was to develop a theory parallel to Galois Theory, 
in which to associate to a differential equation a group and to answer
questions about the dynamic system by studying the correspondent group.

The dynamical system leads to considering a Lie algebra, and the 
solutions correspond to a Lie group.
Citing from \cite{VBLT}, p.601, ``The basic object mediating
between  Lie groups and Lie algebras is the one-parameter
group. Just as an abstract group is a coherent system of
cyclic groups, a Lie group is a (very) coherent system of
one-parameter groups.''.

The first main theorem of LT, {\bf Lie's Exponential Map Theorem} \cite{ELT}, p.82,
states the existence of the Lie algebra of infinitesimal transformations
associated to a given Lie group of transformations
(the details are not essential here). 
This provides a functor
\[
\text{Lie Groups}\rightsquigarrow \text{Lie Algebras}.
\] 
Moreover, if $\mathfrak g$ denotes the Lie algebra of a Lie
group $G$, then there is a natural identification
${\mathfrak g}=T_eG$, and the exponential map $\exp\colon
({\mathfrak g},0)\to (G,e)$ is an isomorphism of germs
of differential manifolds. Then {\bf Lie Main Theorem}
\cite{ELT}, p.83, identifies the Lie algebra structure and
provides the construction of the Lie group given such an
abstract Lie algebra.
More precisely, disregarding issues of
convergence, e.g. considering (pro)nilpotent Lie algebras or
working at a formal level, one can use the exponential map to
\emph{integrate} a Lie algebra ${\mathfrak g}$, i.e., to
construct a Lie group with Lie algebra ${\mathfrak g}$:
\[
{\rm Exp}:\text{Lie Algebras}\rightsquigarrow \text{Lie
Groups}.
\]
Then {\bf Lie Main Theorem} \cite{ELT}, p.83, identifies the Lie algebra structure and
provides the construction of the Lie group given such an abstract Lie algebra.

The important point is the existence of infinitesimal transformations that
generate the finite transformations of the group \cite{ELT}, p.80.
In fact ``For Lie, intuitively, {\em this} is what constitute the continuity
of $G$.'' (loc. cit.; the author's emphasis). 
As Lie wrote, ``A group is called {\em continuous} when all its transformations
are generated by repeating infinitesimal transformations infinitely often ...''.
Sophus Lie had a natural ``proclivity for geometrically flavored intuitive thinking'' (loc. cit.),
consistent with his collaboration with Felix Klein \cite{ELT}, p.3.
Although LT took later an analytic flavor, influenced by the ``analysts and arithmeticians''
from that period, ``topological considerations remained outside Lie's theory
until the mid-1920s when Weyl ... began to introduce them,...'' (loc. cit., p.80).
The return to the geometric flavor of the two ``synthesists'', was achieved by
Cartan and then Chen, as it will be briefly mentioned next.

\begin{rem}\label{R:LTFP}
On the pedagogical side, 
a ``simplified'' account of LT can be achieved by 
{\em not} considering abstract Lie groups,
but only groups concretely realized as groups of matrices \cite{VBLT}, p.601.
On the application side, we interpret groups of matrices 
as representations of quivers with additional structure,
providing the bridge to Feynman Process as representations of more 
general geometric categories: Feynman Categories.
\end{rem}

\subsection{Further developments}\label{SS:LTFD}
The notable developments of LT, which in our opinion leads to DT,
were achieved by Cartan in the late 19th century, S. Bochner in 1946 (formal group laws),
Chen (formal connections; see Remark \ref{R:Chen}), 
Kodaira-Spencer Theory etc. \cite{wiki:CC,wiki:FG}.

{\em Cartan connections} describe the geometry of manifolds modeled on 
{\em homogeneous spaces}, for example Riemann surfaces via the Uniformization Theorem \cite{wiki:UT}.
Recall that Klein's programme suggested that ``geometry'' is the study of a homogeneous space.
``A Cartan geometry modeled on a homogeneous space can be viewed as a {\em deformation}
of this geometry which allows the presence of {\em curvature}'' \cite{wiki:CC}, p.5.

Formal groups mark the transition from topological aspects to an emphasis on geometric aspects.
A {\em formal group} law is a formal power series $F(x,y)$ with coefficients in a ring $R$,
behaving like the product of a Lie group
\cite{wiki:FG}:
$$F(x,y)=x+y+higher\ degree\ terms,$$
e.g. the {\em additive formal group law} $F(x,y)=x+y$, 
Hausdorff series and star products.
Note that the Hausdorff series (Backer-Campbell-Hausdorff formula),
is the source of most techniques achieving deformation quantization,
as proved by \cite{Res}, and underlying Rota-Baxter algebra considerations
leading to Birkhoff decomposition and renormalization \cite{EK}.

The {\em logarithm} of a commutative formal group law is an isomorphism
$f$ between the additive formal group to
$F$:
$$f(F(x,y))=f(x)+f(y).$$
For example the $\chi$ map of \cite{EK} is such a logarithm
(to be explained elsewhere).

The interpretation of formal series as {\em connections} was studied by Chen 
(see \cite{Hueb} and references therein).

The theory of pseudo-groups developed in the early 1900 by Elie Cartan,
was reformulated by Shiing-Shen Chern around 1950s, 
and a general deformation theory for pseudogroups was 
later developed by Kunihiko Kodaira and D. C. Spencer \cite{wiki:PG}.
Further developments and applications of DT includes the proof by Kontsevich of Deligne
conjecture, the Formality Theorem; as proved in \cite{FI},
Kontsevich constructs a non-trivial L-infinity morphism allowing to transfer 
a rather trivial solution of MCE, the Poisson bracket, to another DGLA,
obtaining the desired star-product.

Returning to LT, an essential improvement is the
Milnor-Moore Theorem with the functorial correspondence
between Lie algebras and universal enveloping algebras
(``Lie-Hopf correspondence'').

The next step is to generalize the framework to Rota-Baxter algebras
(R-matrices etc.).

Finally, our next goal following the present article \cite{I-WIFPI}, 
is to explain in detail how a {\em Feynman Path Integral} 
is a further generalization of the multi-valued logarithm (a path integral, really),
making the connection with the theory of {\em Quantum Information}.

\subsection{Lie Theory and Deformation Theory: a preview}
Given a Lie group G and its Lie algebra
${\mathfrak g}= T_eG$, Lie Theory constructs  two maps: 1)
the exponential
$exp: T_eG\to G$ and, in the connected and simply connected case
$\log: G \to T_eG$ provides and inverse.

Now Maurer-Cartan Equation is a substitute for the group $G$ \cite{Getzler},
so the exponential maps can be rewritten as
$exp:T(MC_{\mathfrak g}) \to MC_{\mathfrak g}$,  and the
logarithm goes in the other direction
\footnote{This is actually much more than a bare analogy, see \cite{Getzler}.}.
Written this way, these maps generalize to DGLAs, 
and with some technical effort to $L_\infty$ algebras. 

Namely, let
${\mathfrak g}$ be a  DGLA or, more generally, an
$L_\infty$ algebra, and let $MC_{\mathfrak g}$ be the
Maurer-Cartan equation for ${\mathfrak g}$.  The
tangent space to
$MC_{\mathfrak g}$ is $Z^1({\mathfrak g})$, the space of 
cocycles of ${\mathfrak g}$, and the analog of the 
exponential is a map
$\xy\xymatrix{Z^1({\mathfrak g})
\ar@{-->}[r]& MC_{\mathfrak g}}\endxy$, 
which ``completes'' 
a cocycle with vanishing obstruction to a solution of the
Maurer-Cartan equation; this should be seen as the analog of
the 1-parameter group of Lie Theory.

%
%
\section{Deformation theory and the Huebschmann-Stasheff universal solution}\label{S:dtusol}
The interpretation of the Huebschmann-Stasheff construction of the
universal twisting cocycle solving the Maurer-Cartan equation (MCE) of a DGLA,
is given in terms of the Kuranishi functor.
Comparison with Lie Theory suggests that DT is a higher order version of Lie Theory.
In view of the interpretation by Chen of formal series as {\em connections} \cite{Hueb},
DT is a ``connection geometry'', 
i.e. a graded version of the theory of {\em Cartan connections} \cite{wiki:CC} 
in the context of formal groups,
as it will be explained in \S\S \ref{SS:LT}.
Such a comprehensive study of DT is beyond the aims of the present article.

But first, a brief recall on deformation theory is in order
(\cite{Man}; see also \cite{AIS}).
The main points of Lie Theory to bear in mind will be subsequently reviewed next, 
together with its developments by Chen and Kodaira.
The following subsections contain the technical details, towards a justification of 
the central role of the Kuranishi functor as a linearization device: 
a logarithm for the universal solution of MCE, an analog of a Lie exponential.

\subsection{Recall on deformation theory}\label{S:dt}
The usual presentation of deformation theory focuses on
{\em deformation functors} associated to a DGLA
${\mathfrak g}=\oplus {\mathfrak g}^i$ \cite{Man} p.14, as a
functor of the coefficient (Artinian) ring. We fix such a
``local model'', i.e. a standard model of a pointed formal Lie manifold, 
say ${\mathfrak m}$ the maximal ideal of the
pro-Artinian ring ${\mathbb K}[[h]]$,  and focus on the
functoriality with respect to the  DGLA argument as in Lie theory.

The three main functors are \cite{Man}: 1) the {\em Gauge functor}
$G({\mathfrak g})=exp({\mathfrak g}^0\otimes m)\in Group$, 
2) the {\em Maurer-Cartan functor} of flat connection forms:
$$MC({\mathfrak g})=\{x\in {\mathfrak g}^1\otimes {\mathfrak
m} | dx+\frac12 [x,x]=0\}.$$ To define the third functor, 
note that ${\mathfrak g}\otimes {\mathfrak m}$ is a DGLA,
with
$({\mathfrak g}\otimes {\mathfrak m})^0$
defining an m-adic topology compatible with the algebraic structure 
(h-adic topology of formal power series),
and there is an action of the group
$G({\mathfrak g})$ on $MC({\mathfrak g})$. The corresponding
moduli space is called  the {\em Deformation functor}
$Def({\mathfrak g})=MC({\mathfrak g})/G({\mathfrak g})$.
Following Grothendieck one should study $G({\mathfrak g})$
acting on
$MC({\mathfrak g})$ as a groupoid (more specifically as a
discrete bundle/ local system / flat connection),  with
$Def({\mathfrak g})$ the base space.

Now ${\mathfrak g}^0$ may be identified as the tangent space
to $G({\mathfrak g})$ and the cycles $Z^1({\mathfrak g})$ of
${\mathfrak g}$ as the elements of the tangent space to
$MC({\mathfrak g})$
\cite{Man} p.14.

\begin{rem}\label{R:Chen}
According to Chen \cite{Hueb},
$C={\mathfrak g}\otimes {\mathfrak m}$ is the set of {\em
formal connections}
$D=d+A$, with the {\em perturbation} $A$.
Then $T^1={\mathfrak g}^1\otimes {\mathfrak m}$ is the
tangent space and 
$MC({\mathfrak g})$ is the set of {\em flat connections}
\footnote{Local systems over the formal pointed manifold
with observables $m$.}
$$D^2=0\quad \leftrightarrow \quad \Omega(A)=dA+\frac12[A,A]=0,$$
with $Z^1({\mathfrak g})$ its tangent space at the trivial
connection.
\end{rem}
The functor $Def$ lifts to the derived category,
and a quasi-isomorphism between DGLAs ${\mathfrak g}$ and
${\mathfrak g}'$ induces an isomorphism  between the
corresponding moduli spaces $Def({\mathfrak g})$ and
$Def({\mathfrak g}')$.

\subsection{Kuranishi functor}\label{SS:KF}
The Kuranishi functor \cite{Man} plays an important role
in deformation theory, and in particular when solving MCE \cite{HS},
as noticed in \cite{AI-HS}.

The Kuranishi maps allow to represent deformation functors
$Def({\mathfrak g})$ for which $H^1({\mathfrak g})$ is finite
dimensional.

\newcommand{\To}{\longrightarrow}
Chose a direct sum decomposition corresponding to a complement to the 
space of cycles $C^i$ and a complement to the space of boundaries $B^i$
 (\cite{Man}, p.17;  compare with \cite{AI-HS} Example 6, p.17)
$${\mathfrak g}^i=Z^i\oplus C^i, \quad Z^i=B^i\oplus \C{H}^i$$
and let $h:{\mathfrak g}^{i+1}\to {\mathfrak g}^i$ be the linear map
given by the composition
$${\mathfrak g}^{i+1}\overset{\pi_B^{i+1}}{\To}B^{i+1} 
\overset{d^{-1}}{\To}C^i\subseteq {\mathfrak g}^i$$
where $\pi_B^{i+1}$ is the projection of ${\mathfrak g}^{i+1}$ of kernel
$C^{i+1}\oplus \C{H}^{i+1}$.
Then $x\in \C{H}^i$ iff $d^ix=0$ and $h^ix=0$. 

If follows that $d h+h d=Id-\pi_\C{H}$ (\cite{Man}, p.17, with $H:=\pi_\C{H}$), where 
$$\pi_\C{H}^i:{\mathfrak g}^i\overset{\pi^i}{\to}\C{H}^i \overset{\nabla}{\to}{\mathfrak g}^i,$$
is essentially $\pi^i$, the projection on $\C{H}^i$ of kernel $B^i\oplus C^i$,
and $\nabla$ is the canonical inclusion. 
For later use, note that $h$, an {\em almost contraction} in the sense of \cite{AIS}, p.6,
satisfies the identity:
\begin{equation}\label{E:actr}
(dh+hd)\pi_B=\pi_B, \ \forall x\in B\ dh(x)=x,
\end{equation}
i.e. $h$ inverts $d$ on the complement $C$ of the cycle space $Z$.

We have the following result 
interpreting the Kuranishi map/functor of deformation theory
from the Huebschmann-Stasheff homotopy perturbation theory (HPT) point of view,
as explained in \cite{AI-HS} \S5 (except for taking homology representatives,
for compatibility with \cite{Man}).
\begin{th}\label{T:HSK}
With the above notations, $\C{H}$ is a {\em strong deformation retract of ${\mathfrak
g}$}  (SDR) \cite{HS}, isomorphic to its homology $H_\bullet({\mathfrak g})$:
\begin{equation}
( (\C{H},0)  
\begin{array}{c} {\pi} \\ 
\stackrel{\leftrightharpoons}{\nabla }     
\end{array}  
({\mathfrak g},d),h).\label{E:contr} \end{equation}
\end{th}
Our objective now is to use the Picard method interpretation of the 
Huebschmann-Stasheff construction \S\ref{S:Pic},
allowing to interpret the Kuranishi map as a {\em resolvent of the
Maurer-Cartan curvature}:
\begin{equation}\label{E:MCC}
\Omega=d+1/2[,].
\end{equation}
\begin{defin}
The {\em Kuranishi map} $F:T^1({\mathfrak g})\to T^1({\mathfrak
g})$ is the morphism of functors given by:
$$y=F(x)=x+\frac12 h[x,x], \quad x\in T^1({\mathfrak
g})={\mathfrak g}^1\otimes {\mathfrak m}.$$
\end{defin}
The Kuranishi map is an isomorphism of functors \cite{Man} (Lemma 4.2, p.17).
\begin{defin}
The {\em Kuranishi functor} is given by:
$$Kur({\mathfrak g})=\{y\in \C{H}\otimes {\mathfrak m}
\ \text{ such that }\ p_\C{H}[F^{-1}(y),F^{-1}(y)]=0\}.$$
\end{defin}
In other words $Kur({\mathfrak g})$ is the kernel of the
morphism of functors induced by the composition:
$$\C{H}^1\overset{\nabla}{\to} {\mathfrak g}^1\overset{F^{-1}}{\to}{\mathfrak
g}^1\overset{\Omega}{\to} {\mathfrak g} \overset{\pi}{\to}\C{H}^2.$$ 
We have included the term $dx$ together with $[,]$,
in order to emphasize the MC-curvature $\Omega$ (Equation \ref{E:MCC}),
since the term vanishes under the projection 
$p_\C{H}$ to $\C{H}^2$ (compare with \cite{Man}).
As a benefit, the role of the Kuranishi functor is now more transparent (``vertical forms''
with exact curvature):
$$Kur=Ker\ \pi_\C{H} [(F^{-1})^*(\Omega)].$$
The Kuranishi functor yields a {\em reduction of the MC-equation},
explained next.

Let $N=\oplus N^i$ be defined by:
$$N^i=0, i\le 0, \quad N^1=C^1\oplus \C{H}^1, \quad
N^i={\mathfrak g}^i, i\ge 2.$$ It may be thought of as a
reduction of the space of connections subject to a gauge
condition: 
$$N\cong {\mathfrak g}/B^1, \quad B^1=\{ A\in {\mathfrak
g}^1\ \text{ such that }\ A=df\}.$$
\begin{prop}
The isomorphism $F$ induces an isomorphism \cite{Man}:
\begin{equation}\label{E:KF}
F:MC({\mathfrak g})\cap (N^1\otimes {\mathfrak m})\to
Kur({\mathfrak g}).
\end{equation}
\end{prop}
In other words, under the Kuranishi map $F$,
the Kuranishi functor $Kur$ is isomorphic with
the reduction $MC(N)$.

Now note that the projection ${\mathfrak g}\to {\mathfrak
g}/B^1\cong N$ is a quasi-isomorphism so that the
corresponding deformation functors are
isomorphic:
$$H({\mathfrak g})\cong H(N) \quad \Rightarrow\quad
MC({\mathfrak g})/G({\mathfrak g})\cong MC(N)/G(N).$$ In
conclusion,
\begin{th}\label{T:KF}
For any DGLA ${\mathfrak g}$ the deformation moduli space morphism is \'etale \cite{Man}:
$$Kur({\mathfrak g})\overset{F^{-1}}{\to} MC({\mathfrak g})
\overset{Proj}{\to}Def({\mathfrak g}).$$ 
In particular the Kuranishi functor is locally isomorphic to the  deformation
functor.
\end{th}
Now we can interprete the construction of the universal 
twisting cocycle from \cite{HS}.

\subsection{The Huebschmann-Stasheff construction}
The Huebschmann-Stasheff construction is based on a 
representable version of the Kuranishi functor explained above,
obtained by applying $Hom^{\mathfrak g}=Hom(\cdot,
{\mathfrak g})$. 
In other words, the universal twisting cocycle
$\tau$ solves the MC-equation (master equation) in the 
Chevalley-Eilenberg DGLA associated to the DGLA ${\mathfrak g}$:
$$\diagram
(({\mathfrak g},d),h) \ar@/^2pt/[d]^{\pi} & \lto^{}
S^c(s{\mathfrak g})
\lto_{pr} \dto^{S^c(\pi)} \\ (\C{H}\cong H({\mathfrak g}),0)
\ar@/^2pt/[u]^{\nabla} & (S^c(sH({\mathfrak g})), D)
\ulto^{\tau}\lto^{\tau_{H(g)}},
\enddiagram$$
where $S^c$ is the graded symmetric coalgebra functor 
and $s$ denotes the suspension functor \cite{HS} (see also \cite {AI-HS}).

Recall that $\tau=\tau^1+\tau^2....$ is a universal solution of the MCE in the following sense
(again using representatives from $\C{H}$, instead of homology classes as in \cite{HS,AI-HS}):
\begin{equation}\label{E:UMCE}
\tau^1=\tau_{H({\mathfrak g})}, \quad \Omega(\tau)=0,
\end{equation}
where $\tau^1$ is the ``universal initial condition''.

Explicitly, for any representative $x\in \C{H}$ of a homology class $\pi(x)\in H(\g)$, 
$\tau_x=\tau(x)$ is the solution of the {\em Maurer-Cartan Initial Value Problem} (MC-IVP)
\begin{equation}\label{E:MCIVP}
d\tau_x+\frac12 [\tau_x, \tau_x]=0, \quad \tau_x^1=x.
\end{equation}
Now the relation satisfied by the contraction $h$ (Equation \ref{E:actr})
together with the fact that $\tau_x^1=x\in Z$ while $\tau_x^b\in C^b, b>1$,
imply that the MC-IVP \ref{E:MCIVP} is equivalent to 
(see also \cite{AI-HS}, p.17, modulo an irrelevant change of sign in MCE):
\begin{equation}\label{E:MCEC}
\tau_x+\frac12 h[\tau_x,\tau_x]=x.
\end{equation}
Now define the {\em path integral starting at $x$}:
$$C_x(y):=x-h[y,y].$$
Since $h[x,x]$ is an {\em h-adic topology contractions}, so is $C_x$.
Then Equation \ref{E:MCEC} becomes:
\begin{equation}\label{E:MCFP}
C_x(\tau_x)=\tau_x.
\end{equation}
To better understand the role of the universal solution,
as the analog of the exponential of LT,
and that of the Kuranishi map (natural transformation) as the analog of the logarithm,
we will review the well known method of Picard for solving
initial value problems.

\subsection{Picard's method and almost contractions}\label{S:Pic}
The use of such contractions \cite{HS} 
was ``rediscovered'' by the present author
in \cite{AIS}, and called ``almost contractions'' in contrast with
{\em contracting homotopies} \cite{HSta}, p.125.

We will recall the idea behind Picard's method \cite{Kr} p.285
and try to recast it in terms of
homotopical algebra.

The initial value problem (IVP) is equivalent to the integral solution:
$$dy/dt=K(y,t), y(0)=y_1\quad y=y_1+\int_0^x K(y(t),t)dt,$$
which can be found iteratively:
$$y_{n+1}=y_1+\int_0^x K(y_n,t)dt, n\ge 1,$$
provided $K$ is Lipschitz (contraction)
\footnote{This is no longer a constraint when looking for formal solutions.}.

In the above example, $d/dt$ and $\int_0^x$ are reminiscent of 
such a pair for the algebra of formal power series $k[[x]]=\oplus P^n$:
$$d(x^n)=nx^{n-1}, h(x^n)=\int_0^x t^n dt= x^{n+1}/(n+1).$$
Then $d h= Id$ and $h d=Id - H$, 
where the evaluation at 0 $H(f)=1/2 f(0)$, 
is the projection on the subspace of constants, 
of kernel $m$ (augmentation).

Then $dh+hd=2Id-H$, except of course $d^2\ne 0$ etc., and 
Picard's solution is now:
$$y_{n+1}=y_1+h[F(y_n,t)].$$

\subsection{The exponential and logarithm of a DGLA}
The analog for formal deformations consists of Maurer-Cartan equation:
$$d\tau=K(\tau)=1/2[\tau,\tau], \tau(0)=\tau^1$$
with a similar iterative solution:
$$\tau_{n+1}=\tau^1+h[K(\tau_n)], $$
Let us first briefly recall the construction from \cite{HS}, p.10,
applying the Theorem 2.7 (loc. cit. p.9) to our contraction (SDR)
from Equation \ref{E:contr}; 
for further details see \cite{AI-HS}, p.21.

As mentioned above (see also \cite{HS}, Addendum 2.8.1, p.10),
the twisting cochain is determined by
$\tau:S_D^c[sH({\mathfrak g})]\to {\mathfrak g}$,
as an element of degree -1 of the Chevalley-Eilenberg
DGLA $Hom(S^c_D[sH({\mathfrak g})],{\mathfrak g})$, which
satisfies the master equation 
$$D\tau=1/2[\tau,\tau],$$
with $D$ the coderivation extending $d$.
The solution is constructed iteratively:
$$\tau=\tau^1+\tau^2+...+\tau^b+..., \quad
\tau^j:S^c_j{\mathfrak g}\to {\mathfrak g}, j\ge 1$$
$$\tau^{b}=\frac12 h\sum_{i=1}^{b-1}[\tau^i,\tau^{b-i}].$$
Consider now the {\em h-adic contraction}:
$$C(x)=x-\frac12 h[x,x],$$
together with a sequence defined recursively by:
\begin{equation}\label{E:yn}
y_1=\nabla, \quad y_{n}=C(y_{n-1}), n\ge 2.
\end{equation}
Then we obtain the following interpretation of Stasheff-Huebschmann construction
together with the Kuranishi map as its inverse.
\begin{th}\label{T:log}
The universal Maurer-Cartan solution $\tau$ is the fixed point $\tau=\lim y_b$
of the h-adic topology contraction $C$ associated to the contraction $h$
and the initial value $\nabla:\C{H}\to {\mathfrak g}$.

Explicitly, 
the solution $\tau_x$ of the Maurer-Cartan initial value problem Equation \ref{E:MCIVP} is 
the fixed point of the h-adic contraction $C_x$ determined by the
contraction $h$ and the initial value $x$.

Moreover the Kuranishi map $F(x)=2I(x)-C_x(x)$ is its inverse,
recovering the initial value of the solution of the MC-IVP:
$$F(\tau_x)=x.$$
\end{th}
\begin{pf}
First note that $C$ is a contraction relative to the 
grading ($h$-adic) topology of $S^c(sH({\mathfrak g}))$.
Now the $b$-th component of the fixed point $\tau$ 
satisfies the above recursive relation:
$$\tau^b=C(\tau)^b=\frac12(h[\tau,\tau])^b=
\frac12 h(\sum_{i=1}^{b-1}[\tau^i,\tau^{b-i}]).$$
Now $\tau$ is the h-adic limit of $y_n=\sum_{b=1}^n \tau^b$, satisfying the 
recursive relation \ref{E:yn}.

The h-adic convergence is clear, so we only note that:
$$y_1=x, \ \lim y_{n+1}=C(\lim y_n) \quad \Rightarrow \quad x=F(\tau_x).$$
\end{pf}
\begin{defin}\label{D:explog}
The universal solution $\tau$ is called the {\em exponential of the DGLA} $\g$,
and denoted $exp_\g$, 
while the Kuranishi map $F$ is called the {\em logarithm of the DGLA} $\g$,
denoted $\log_\g$
\footnote{See \S\ref{S:conclusions} for additional motivation for the terminology.}.
\end{defin}
\begin{rem}
In view of Theorem \ref{T:KF}, 
we think of the contraction $C$ 
associated to the splitting homotopy $h$,
as the projection on the moduli space $Def_L=MC_L/G_L$.
After ``gauge fixing'', i.e. removing the redundant boundaries
by restricting to $C^1\oplus \C{H}^1$ (see \cite{Man}, p.18),
the Kuranishi functor $F=Id\pm C$ becomes an isomorphism (see Equation \ref{E:KF}).
\end{rem}
\begin{rem}
The use of an almost contraction in \cite{AIS} when solving 
the MCE $[*,*]$ for an associative star-product $*$ is 
prototypical of the recursive construction of Stasheff-Huebschmann
in view of the fact that any DGLA $(L,d,[,])$
can be augmented to a {\em pointed DGLA} $(L_d, ad_d,[,]_d)$, 
by adjoining the derivation $d$ \cite{Man}, p.5.
Then MCE $dx+\frac12[x,x]=0$ is equivalent to $[x,x]_d=0$,
which can be solved as in \cite{AIS} using an
almost contraction,
or as in \cite{HS}.
\end{rem}
We have interpreted Equations (2.7.2) and (2.7.3) from \cite{HS}
as giving a fixed point of a contraction, in order to
to emphasize the typical approach of solving a differential (or integral) 
equation iteratively, as it will be recalled next.

\subsection{Higher Dimensional Lie Theory}\label{S:hlt}
The idea is that {\em Deformation Theory is a Higher Lie Theory}.
In Lie theory the infinitesimal Lie algebra is exponentiated to obtain
the closed 1-parameter groups $a_t=E^{At}$. 
In the classical Lie groups case,
the exponential is the solution of the initial value problem
$$dy/dt=Ay, \quad y(0)=e.$$
The solution can be obtained as a formal power series $y=\sum a_n x^n$,
which leads to recursive formulas for analytic coefficients
$a_{n+1}=f(a_n)$, or using Picard's method by solving the equivalent 
integral equation $y=y_0+\int Ay$, which also leads to an iterative procedure:
$$y_{n+1}=y_0+\int A y_n.$$
Now $\int$ plays the role of the contraction for $d$: 
$$d\int-\int d=Id-H,$$
so the Huebschmann-Stasheff construction may be thought of as 
{\em the Picard's method for finding the exponential of a higher Lie theory}.

\subsection{Generalization to the L-infinity case}\label{S:linfcase}
Now $K(y)=[y,y]$ is a contraction in the h-adic completion sense,
so the (formal) series converges.

The construction of the solution can be generalized to the full MC-equation,
corresponding to an $L_\infty$-algebra $(g,d,[,],[,,],...)$:
$$d\tau+\sum_{n\ge 2} [\tau,...,\tau]/n!=0,\qquad (Q(\tau)=0).$$
We interpret this equation as a higher version of Lie theory,
which corresponds to $dy=0, \quad d=D+A$.

Then the solutions of MCE correspond to a pointed formal manifold,
generalizing the case of a DGLA, with its exponential and logarithm.
In this sense 
{\em deformation theory is a higher version of Lie theory}.

\section{The Laplacian: doubling or halving}\label{S:dorh}
The homotopy $h$ used to construct a solution of MCE
is determined by the splitting of $\g$ into homology $\C{H}$,
boundaries $B$ and residual piece $C$, 
reminiscent of the Hodge decomposition of the de Rham complex
of a Riemannian manifold.

\subsection{The *-operator of a contraction}\label{S:hd}
In \cite{AI-HS}, p.16. it is proved that 
in general a SDR $N\to M$ with contraction $h$ defines a *-operator
$*=h+d_N+Id_{\nabla}$, such that the codifferential 
$$d^*=*d*^{-1}=h$$ 
is the contraction $h$ and the associated Laplacian 
$$\Delta=(d^*+d)^2=\nabla\pi-Id_N$$
is essentially a projection.

In our context, with $\g\to H(\g)$ a SDR with contraction $h$ corresponding to 
the decomposition
$$\g=\C{H}\oplus B \oplus C,$$
the ``Hodge'' isomorphism is 
$$*=(h+d)\oplus Id_\C{H},\quad d^*=h,$$
and the ``Hamiltonian''
$${\bf H}=-\Delta=Id_\g-\nabla\pi=Id_\g-H$$
is the projection onto the space of boundaries and {\em coboundaries}
$$\C{D}(B)=B\oplus B^*, \quad B^*=Im\ d^*= Im\ h=C,$$
of kernel $\C{H}$: the {\em harmonic forms}.

\begin{defin}\label{D:hd}
The {\em Hodge decomposition} associated to a contraction $h$ of a SDR data is:
$$\g=B\oplus \C{H}\oplus B^*, \quad B=Im d,\ B^*=Im d^*,$$
where $d^*=*d*=h$ and $*$ is the associated Hodge operator.
\end{defin}
As expected, the cycles are $Z=ker d=B\oplus \C{H}$ and cocycles
$Z^*=ker d^*=B^*\oplus \C{H}$.
\begin{rem}\label{R:transfer}
It is interesting to note that $g$ has the structure from
Lemma 1 of \cite{Kon2}, p.12.
$$\g\cong H(\g)\oplus B\oplus B^*.$$
Recall that the basic principle of homotopy perturbation theory,
the {\em Gugenheim principle},
refers to the transfer of structure under quasi-isomorphisms (\cite{Lambe}, p.2-3):
given a resolution $\epsilon:A\to H(A)$ (quasi-isomorphism)
and a deformation of $H(A)$, is there a deformation of $A$
transfered via the quasi-isomorphism?
$$\diagram 
P-algebra\ A \rto^{\quad \epsilon} \ar@{|.>}[d] & H(A) \ar@{|->}[d]^{Deformation} \\
A[[h]] \rto^{\epsilon[[h]]} & H(A)[[h]].
\enddiagram$$
Kontsevich answers affirmatively this question (loc. cit.) 
in the case of $P$-algebras for certain operads $P$ 
satisfying some technical conditions.
\end{rem}
Returning to the Hodge decomposition, 
$$**=Id_\C{H}-Id_{\C{D}(B)},$$
is a sort of a {\em Hilbert transform}.
On $\C{D}(B)$, i.e. outside the non-trivial piece $H(g)$,
the star operator $*$ is also a complex structure
or a ``Dirac structure''
\footnote{Generalized complex structures are a 
complex analog of Dirac structures \cite{GCS}.}
$$**=-Id_{\C{D}(B)}=\Delta.$$
Of course one could take $-h$ as a contraction instead,
and obtain the Laplacian as the projection (positive operator) 
and investigate the relation with the harmonic oscillator $[d,d^*]=Id_{\C{D}(B)}$.
\begin{rem}
We think of a decomposition of $\g$ as above as a 2-charts 
bundle atlas; it is obtained by doubling $B$ via 
the complex/Dirac structure $*$ 
and then gluing the two pieces along $\C{H}$:
$$\diagram
0 \rto & H(\g) \ar@/^2pt/[r]^{\nabla} & \ar@/^2pt/[l]^{\pi} \g \ar@/^2pt/[r]^{-\Delta} & 
(\C{D}(B),*) \ar@/^2pt/[l]^{} \rto & 0.
\enddiagram$$
\end{rem}
The deformation of the initial DGLA structure on $\g$ to an $L_\infty$-algebra structure
is a ``trivialization of the bundle'': 
it splits the extension in the larger homotopy category.
\begin{th}\label{T:splitting}
The contraction $h$ of a SDR $\g\to H(G)$ 
determines a Hodge decomposition $(\g,d,d^*=h, *=d+h)$ of the DGLA $g$.
Relative to the corresponding $L_\infty$-algebra structure on the homology
obtained by deformation,
$h$ deforms to a splitting contraction of $L_\infty$-algebras:
$$\g\cong \C{D}(B)\oplus \C{H}.$$
\end{th}
Alternatively, the SDR data for the (DG) Lie algebra $\g$ 
prompts to interpret the quasi-isomorphism $\g\overset{\epsilon}{\to} H(\g)$
as an augmentation with augmentation ideal 
the ``double of boundaries'' $(\C{D},*)$ (see Remark \ref{R:transfer}):
$$\diagram
& 0 \rto \dto & Z \drto^{[\ ]} \ar@{^(->}[d] & & \\
0 \rto &  (\C{D}(B),*) \ar@{^(->}@<2pt>[r] & 
\g \ar@<2pt>[l]^{\quad -\Delta} \ar@<2pt>[r]^{\pi} & 
H(\g) \ar@<2pt>[l]^{\nabla} \rto &  0.
\enddiagram$$
Then the restriction of the Hodge operator $*$ (Dirac structure) 
to this double is a complex structure:
$$J=*_{|\C{D}(B)}, \quad *^2=\Delta.$$
The relation with quantization is considered next.

\subsection{Gauge reduction and the Hodge-Cartan decomposition}
The moduli space $Def_\g$ of deformations modulo equivalences (gauge transformations)
is the ``Hamiltonian reduction'' relative to the Hodge decomposition
associated to a contraction satisfying the side conditions.
The integrability condition, i.e. the absence of obstructions to deformations,
is a Cartan condition for the decomposition of the Lie algebra of cycles into 
harmonic forms and exact forms.

Essentially $d^*=*d*=h$, and therefore $C=B^*$ is the space of cocycles
and the Hodge decomposition of the DGLA $\g$ determined by the contraction $h$ (SDR) is:
$$\g=B\oplus \C{H} \oplus B^*.$$
Note that the absence of obstructions can be traced to the fact that the 
Hodge decomposition is also a ``Cartan decomposition condition'':
$$d^*[\C{H},\C{H}]\subset B^*.$$
The relation with the $dd^*$-Lemma will be investigated elsewhere.

\subsection{Remarks on deformation quantization}\label{S:DQ}
The above interpretation of the deformation of the initial structure
to allow a splitting suggests a connection with deformation quantization.

One approach, the so called {\em bialgebra (deformation) quantization}
uses the additional bialgebra structure for twisting and producing
a quantization (quantum groups as Hopf algebra deformations etc.
\cite{CP}).
The point that deformation quantization via Hopf algebra deformations 
underlies the renormalization process in the algebraic framework of
Connes-Kreimer was already made in \cite{IM}.
As briefly mentioned in \cite{I-PROPF,I-FL},
Feynman rules and renormalization aims to represent a Feynman category
(DG-coalgebra PROP) as a ``quotient'' of the universal PROP generated by the
(2-pointed) Riemann sphere $\B{C}P^1$, 
with underlying algebraic-geometry object we call ``bifield'' (compare with
Hopf/bialgebra at the infinitesimal level).
With this in mind, deformation theory is a ``higher Lie theory''
targeting Lie bialgebras (Hopf algebras) and the corresponding
algebraic-geometric picture.
We think that the ``bifield'' plays the role of a quantum information propagator,
and the involution $z\to 1/z$, which maps ``sequential addition'' into the 
``parallel addition'' of inverses, is fundamental in the quantum computing
interpretation of space-time \cite{I-ST}.

\section{Conclusions and further developments}\label{S:conclusions}
It was explained that deformation theory exponentiates solutions
of Maurer-Cartan solutions corresponding to deformations of algebraic structures
in a similar manner with Lie theory. 
The relevant correspondences are given by the Kuranishi map and functors.
They are related to the universal solutions of Huebschmann and Stasheff \cite{HS},
which are universal twisting cocycles.

The role of such a deformation in the context of a strong deformation retract
was explained: the ``almost contraction'' of \cite{AIS}
or the contraction $h$ of \cite{HS} provide a ``Hodge duality''
which splits in the larger category of deformed structures (L-infinity algebras).

Specifically, the universal solution $\tau$ is the analog of the exponential from LT.
Its value $\tau(x)$ on a cocycle $x\in Z=T^1(MC)$ 
is the unique solution $\tau_x$ of the MC-IVP; it is a formal 1-parameter deformation 
in the direction of $x$.
The inverse of the exponential $\tau_x\to x$ is given by the Kuranishi map
$F(\tau_x)=x$, so the Kuranishi functor $Kur(\g)$
is the analog of the Lie algebra functor, with value the tangent space at zero
to the non-linear space of solutions of the MCE.
Moreover the Hodge-Cartan decomposition of the Lie algebra $\g=B^*\oplus \C{H}\oplus B$,
associated to the contraction $h$
induces a ``Hamiltonian'' reduction of the moduli space of deformations $Def=MC/G$
modulo gauge transformations $G$.
\begin{th}\label{T:main}
The Huebschmann-Stasheff universal solution of the MCE determined by 
a SDR-data is the inverse of Kuranishi map.
$$\diagram
Z(\g)=T^1(MC)\ar@{->>}[d]_{\pi_\C{H}} \rto^{\quad exp:=\tau} & 
MC(\g)\ar@{->>}[d]^{/G} \\
\C{H}\cong Kur(\g) & \lto^{\quad \log:=F} Def(\g).
\enddiagram$$
The Kuranishi functor is the linearization of the deformation functor,
compatible with the ``Hamiltonian reduction'' of the MCE 
corresponding to the Hodge-Cartan decomposition of the SDR-data.
\end{th}
More general, the Hodge decomposition $\g=B\oplus\C{H}\oplus B^*$ 
corresponding to the homotopy $d^*$, 
appears in the context of the $dd^*$-lemma \cite{dd}:
$$dd^*+d^*d=Id - (Id-\Delta),$$
where $\Delta$ is the corresponding Laplacian,
with applications to the deformation theory of complex structures.

Another interesting application regards Dirac structures and 
generalized complex structures \cite{GCS},
which include the symplectic and the complex case in a common framework,
with possible implications to mirror symmetry.
The above context of a SDR data for $g$ allows for an additional generalization 
from generalized complex structures on $T\oplus T^*$ 
to the Riemann-Hilbert problem framework, 
with applications to renormalization,
as hinted in \S \ref{S:DQ}.

Returning to the above ``doubling and gluing'' process,
the relation with bialgebra quantization emerges (\S \ref{S:DQ}).
The splitting of the $g$ after deformation should exhibit a Lie (L-infinity)
bialgebra structure with an r-matrix corresponding to the involution $*$.
The ``universal object'' in this context is the Riemann sphere bifield $\B{C}P^1$
with its Hopf algebra of functions $A=Hom(\B{C}P^1,\B{C})$ and universal
non-commutative de Rham complex (DG-Hopf algebra) $\Omega^\bullet(A)$.
The relation between a Hodge structure $(d,^*)$, $dd^*$-lemma or SDR 
and the bialgebra structure from the perspective 
of the above universal non-commutative Hodge-de Rham complex,
will be investigated elsewhere.

Regarding the overall picture,
i.e. deformation theory with its applications to quantization,
it is a natural continuation of Lie Theory and universal 
enveloping algebras, which in turn leads beyond Hopf algebras,
to Rota-Baxter algebras and its applications to renormalization \cite{I-RBA}.

\section{Acknowledgments}
I would like to thank Fusun Akman for useful discussions and to Domenico Fiorenza 
for comments which improved a critical part of the article.


\label{lastpage}

\end{document}